\documentclass[preprint]{elsarticle}
\usepackage{amsmath, amsthm, amssymb}
\usepackage{color}
\usepackage{subcaption}
\usepackage{verbatim}
\usepackage[all]{xy}
\newcommand{\eref}[1] {Equation~(\ref{eqn:#1})}
\newcommand{\fref}[1] {Figure~\ref{fig:#1}}
\newcommand{\tref}[1] {Table~\ref{tbl:#1}}
\newcommand{\dref}[1] {Def.~\ref{def:#1}}
\newcommand{\pref}[1] {Prop.~\ref{pro:#1}}
\newcommand{\yref}[1] {Hyp.~\ref{hyp:#1}}

\newtheorem{prop}{Proposition}

\newtheorem{hyp}[prop]{Hypothesis}
\newdefinition{df}[prop]{Definition}
\newdefinition{rmk}[prop]{Remark}
\begin{document}

\begin{frontmatter}

\title{Enumeration of octagonal tilings}

\author[chicago,cmu]{M. Hutchinson\corref{cor1}}
\ead{maxhutch@uchicago.edu.edu}
\address[chicago]{Department of Physics, University of Chicago, 5720 S. Ellis Avenue, Chicago, IL, USA, 60637}
\author[cmu]{M. Widom}
\ead{widom@andrew.cmu.edu}
\address[cmu]{Department of Physics, Carnegie Mellon University, 5000 
Forbes Ave, Pittsburgh, PA, USA, 15213}
\cortext[cor1]{Corresponding author}

\begin{abstract}
Random tilings are interesting as idealizations of atomistic models of quasicrystals and for their connection to problems in combinatorics and algorithms.  
Of particular interest is the tiling entropy density, which measures the relation of the number of distinct tilings to the number of constituent tiles.
Tilings by squares and 45$^{\circ}$ rhombi receive special attention as presumably the simplest model that has not yet been solved exactly in the thermodynamic limit.
However, an exact enumeration formula can be evaluated for tilings in finite regions with fixed boundaries.
We implement this algorithm in an efficient manner, enabling the investigation of larger regions of parameter space than previously were possible.
Our new results appear to yield monotone increasing and decreasing lower and upper bounds on the fixed boundary entropy density that converge toward $S_\infty = 0.36021(3)$. 
\end{abstract}


\end{frontmatter}

\section{Introduction}

A tiling covers a space with a set of compact figures that fill the
space without gaps or overlaps~\cite{Grunbaum1987}.  Often the tiles can
be arranged in many distinct patterns, leading to an ensemble of
tilings with finite entropy density per tile or per area.  Such
tilings are called random tilings, and have been deeply studied within
the physics, mathematics and computer science communities, as we
outline below.

Random tilings entered physics as statistical mechanical models of
dimers~\cite{Kasteleyn1963} and of rough solid surfaces~\cite{Blote1982}.  With the
discovery of quasicrystals~\cite{Shechtman1984,Levine1984}, random tiling
models~\cite{Elser1985,Henley1991,Li1992} were proposed as a natural mechanism to
explain the emergence of quasiperiodicity without relying on
constraints, such as matching rules that define the Penrose tiling.
Because of their high entropy, random tiling models also provide a
simple mechanism to explain thermodynamic stability against competing
crystal phases~\cite{Widom1987}.

Owing to the infinite strengths of their interactions, random tiling
models may have unconventional thermodynamic limits.  Indeed, for
certain tilings inside polygonal fixed boundaries, the entropy density
is a function of the boundary shape~\cite{Elser1984} and is spatially
nonuniform~\cite{Dest98}. The ``arctic circle''
phenomenon~\cite{Propp1996,Cohn2000} provides a striking example in which
the entropy vanishes outside a circle inscribed within the boundary.
Free or periodic boundary conditions generally restore the
conventional thermodynamic limit, and in one noteworthy case periodic
boundary conditions enable an exact solution for the
entropy~\cite{Widom1993,Kalugin1994}.

Fixed polygonal boundaries can be more convenient for practical
numerical calculations in some cases~\cite{Dest04}.  One advantage of
this geometry is a lifting of a $d$-dimensional rhombus tiling into a
higher $D$-dimensional space~\cite{Henley1991,Li1992} with $D>d$.  The
ensemble of tilings corresponds to fluctuations of a $d$-dimensional
directed hyper-surface in the $D$-dimensional hypercubic crystal
lattice.  Each such surface can be placed into 1-1 correspondence with
$d$-dimensional integer partitions.
For these reasons, such problems are called $D\rightarrow d$ problems.

\subsection{Problem definition}

Here, we address the $4 \rightarrow 2$ problem.
This problem can be simply stated: 
given an octagonal region of integral centro-symmetric side lengths $a, b, c, d$ and $N$ tiles consisting of squares and 45 degree rhombi of unit edge length, 
how many configurations, $\Omega$, fill the octagon without overlapping each other or the octagonal boundaries?

\begin{df}[Octagonal tiling]
An octagonal tiling is a non-overlapping covering of a centro-symmetric octagon with integer edge lengths $a,b,c,d$ by the unit square and its $45^\circ$ rotation and the unit $45^\circ$ rhombus and its $45^\circ, 90^\circ$ and $135^\circ$ rotations.
\end{df}
%Henceforth, a space filling, non-overlapping configuration of tiles within a fixed-boundary region is referred to as a `tiling,' and sets of tilings are referred to as ensembles.
%\begin{rmk}
The number of tiles in an octagonal tiling is:
\begin{equation}
N = ac + bd + (a+c)(b+d)
\end{equation}
%\end{rmk}

\begin{df}[$\alpha, \beta$ parameterization]
We specialize to the case of $b = d = \alpha$, $a = c = \beta \alpha$, and $\beta \ge 1$.
In this parameterization, $\alpha$ controls the size of the region and $\beta$ provides some control on the shape.
\end{df}

\begin{df}[Octagonal ensemble]
The octagonal ensemble is the set of all tilings of an octagon.
The cardinality of the ensemble is denoted $\Omega^{oct}(a,b,c,d)$ or just $\Omega(a,b,c,d)$.
Using the $\alpha, \beta$ parameterization, we also write $\Omega(\alpha, \beta) = \Omega(\beta \alpha, \alpha, \beta \alpha, \alpha)$.
\end{df}

\begin{df}[Entropy density] \label{def:S}
The entropy density is defined as:
\begin{equation}
\label{eqn:S}
S(a,b,c,d) = \frac{1}{N}\ln{\Omega(a,b,c,d)}
\end{equation}
Using the $\alpha, \beta$ parameterization, we also write $ S(\alpha, \beta) = S(\beta \alpha, \alpha, \beta \alpha, \alpha) $.
\end{df}

\begin{hyp}[Thermodynamic limit]
The thermodynamic limit is the limit as the number of tiles grows while preserving the ratios of the edge lengths.
Formally, it equals:
$$ \lim_{\alpha \rightarrow \infty} S(\alpha, \beta) = f(\beta) $$
where $f(\beta)$ is an unknown function.
The thermodynamic limit is widely believed to exist\cite{Elser1984,Dest98,Dest01} for all $\beta > 0$.
\end{hyp}

\subsection{Previous work}
There is an exact polynomial-time expression for the cardinality of the ensemble of hexagonal tilings.
\begin{df}[Hexagonal tiling]
A hexagonal tiling is a non-overlapping covering of a centro-symmetric hexagon with integer edge lengths $a,b,c$ by the unit square and the unit $45^\circ$ rhombus and its $45^\circ$ rotation.
\end{df}
%\begin{rmk}
The number of tiles in an hexagonal tiling is:
\begin{equation}
K = ab + ac + bc
\end{equation}
%\end{rmk}
%\begin{rmk}
Definitions of hexagonal tilings are recovered from octagonal tilings when setting $d = 0$.
%\end{rmk}

\begin{prop}[MacMahon Formula] \label{pro:macmahon}
The cardinality of the ensemble of hexagonal tilings of the region $(a,b,c)$ is given by~\cite{MacMahon}:
$$ \Omega^{hex}(a,b,c) = \prod_{i=1}^a \prod_{j=1}^b \prod_{k=1}^c \frac{i + j + k - 1}{i + j + k - 2} $$ 
\end{prop}

Two methods for the enumeration of $4 \rightarrow 2$ tilings have been proposed by Destainville et al.:
an algebraic formulation~\cite{Dest04} that generalizes limited results by Elnitsky~\cite{Elnitsky1997}, 
and recursive formulation~\cite{Dest01} which takes advantage of the hypersurface fluctuation structure.
The time complexity of the algebraic formulation scales as $\exp(\gamma_1 N)$ while the the recursive formulation scales as $\exp(\gamma_2 \sqrt{N})$, for some $\gamma_1, \gamma_2$, 
as we will show, making the recursive formulation more efficient for the enumeration of large tilings.

\subsection{Outline}

First, we review an algorithm originally devised by Destainville {\it et al.}~\cite{Dest01}.
Next, we briefly discuss the efficient implementation of said algorithm, framing it as a breadth-first graph operation and showing it to have favorable complexity compared to a more recent algorithm, also by Desstainville {\it et al.}~\cite{Dest04}.
Then, we present new enumerative results for both the entropy density and a constituent sequence $a_j$.
Finally, we propose relationships within the $a_j$ sequence that provide new bounds on the thermodynamic limit for the entropy density.

\section{Methods}
\begin{figure}
\begin{center} 
\includegraphics[width=0.9\textwidth]{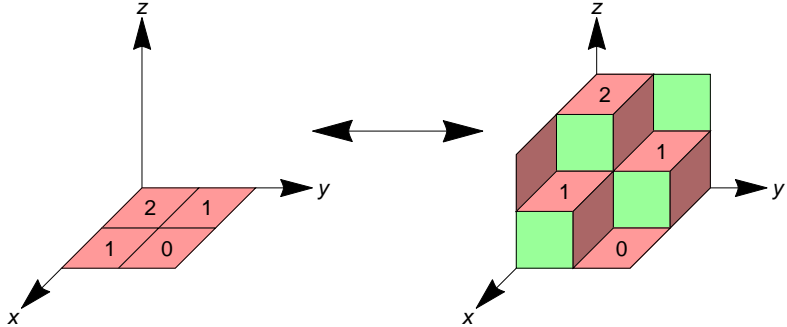}
\end{center}
\caption{ {\bf Example equivalence.} Plane-partitions and hexagonal tilings for size $(2,2,2)$.}
\label{fig:planepartition}
\end{figure}

We begin by demonstrating the construction of a single octagonal tiling of a region with size $a,b,c,d$ as a cascade of partition problems.
The first step is the construction of a plane partition of size $a,c$ and height $b$.

\begin{df}[Plane partition]
A plane partition of size $a,c$ and height $b$ is a matrix $[h_{i,j}]\in \mathbb{N}^{a \times c}$, where
\begin{itemize}
  \item $h_{i,j} < b$
  \item $h_{i,j} \le h_{i-1,j}$
  \item $h_{i,j} \le h_{i,j-1}$
\end{itemize}
\begin{prop}
There is a bijection between plane partitions of size $a,c$ and height $b$ and hexagonal tilings of size $a,b,c$.
\end{prop}
\end{df}
This is equivalent to selecting a membrane from the faces in an $a\times c \times b$ cubic lattice that connects the points $(0,0,b)$ and $(a,c,0)$ 
and uses faces with non-negative normals.
An example for $a=c=b=2$ can be seen in \fref{planepartition}.

We can generalize the notion of a plane partition to generalized partitions.
\begin{df}[Generalized partition]
A generalized partition is a map from a partially ordered set, $P$, to the naturals, $L: P \rightarrow \mathbb{N}$, such that for $x,y\in P$, if $x \succeq y$ then $L(x) \le L(y)$.
The height of $L$ is the greatest element of the image of $L$. 
The set $p_i = \left\{x \in P : L(x) = i\right\}$ is the i-th part. 
\end{df}
%\begin{rmk}
Plane partitions are generalized partitions with $P = \left\{(i,j) \in \mathbb{N}_a \times \mathbb{N}_c \right\}$, with the partial ordering $(i,j) \preceq (k,l)$ iff $i \le k$ and $j \le l$.
%\end{rmk}

We construct a second partition problem on the hexagonal tiling.
\begin{df}[Hexagonal partition]
A hexagonal partition of size $a,b,c$ and height $d$ is a generalized partition on a hexagonal tiling of size $a,b,c$ with height $d$.
The elements of the hexagonal tiling are faces in a cubic lattice.
The partial ordering is the transitive closure of a binary relation between faces that share an edge.
\begin{itemize}
  \item If the edge is parallel to the z-axis, then the face with the greater average $y$ coordinate is less than or equal to the face with the lesser average $y$ coordinate.
  \item If the edge is parallel to the y-axis, then the face with the greater average $z$ coordinate is greater than or equal to the face with the lesser average $z$ coordinate.
  \item If the edge is parallel to the x-axis, then the face with the greater average $y$ coordinate is less than or equal to the face with the lesser average $y$ coordinate.
\end{itemize}
\begin{prop}
There is a bijection between hexagonal partitions of size $a,b,c$ and height less than or equal to $d$ and octagonal tilings of size $a,b,c,d$.
\fref{oct} illustrates this.
\end{prop}
\end{df}

\begin{figure}
\begin{center} 
\includegraphics[width=0.9\textwidth]{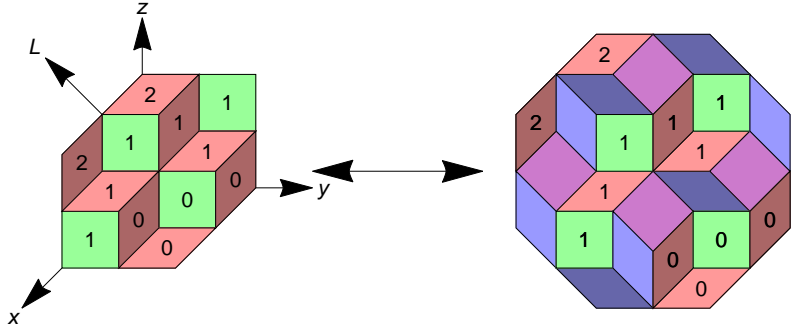}
\end{center}
\caption{ {\bf Octagonal tiling.} Size $(2,2,2,2)$.}
\label{fig:oct}
\end{figure}

\begin{figure}[h]
\centering
\begin{subfigure}[b]{0.45\textwidth}
\includegraphics[width=0.9\textwidth]{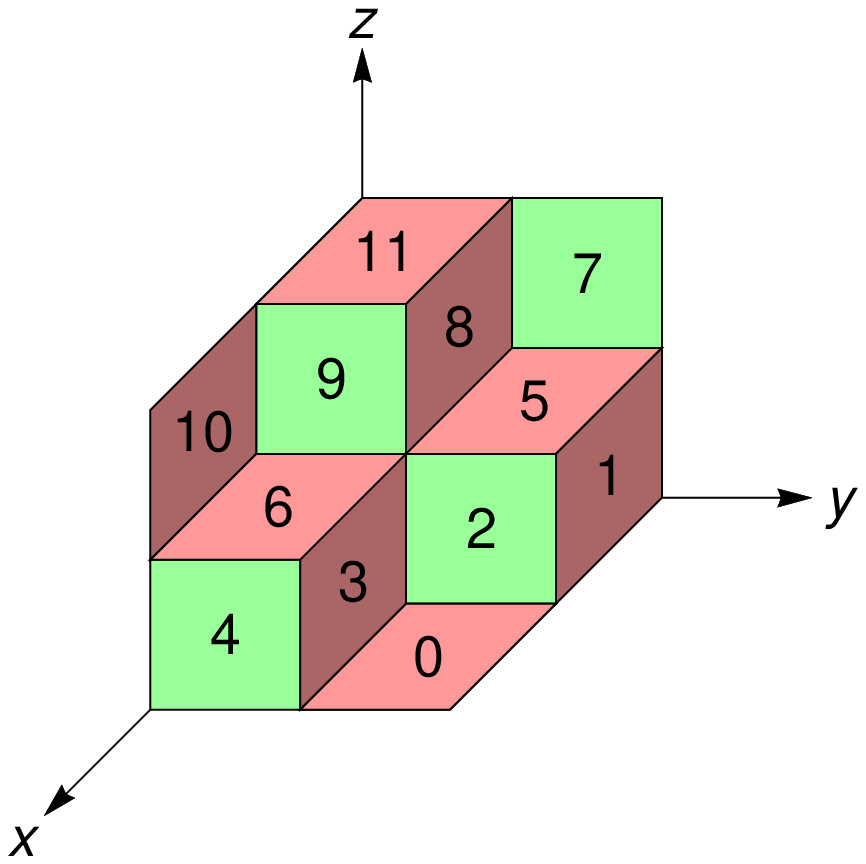}
\caption{\label{fig:ordera}}
\end{subfigure}
\begin{subfigure}[b]{0.45\textwidth}
\includegraphics[width=0.9\textwidth]{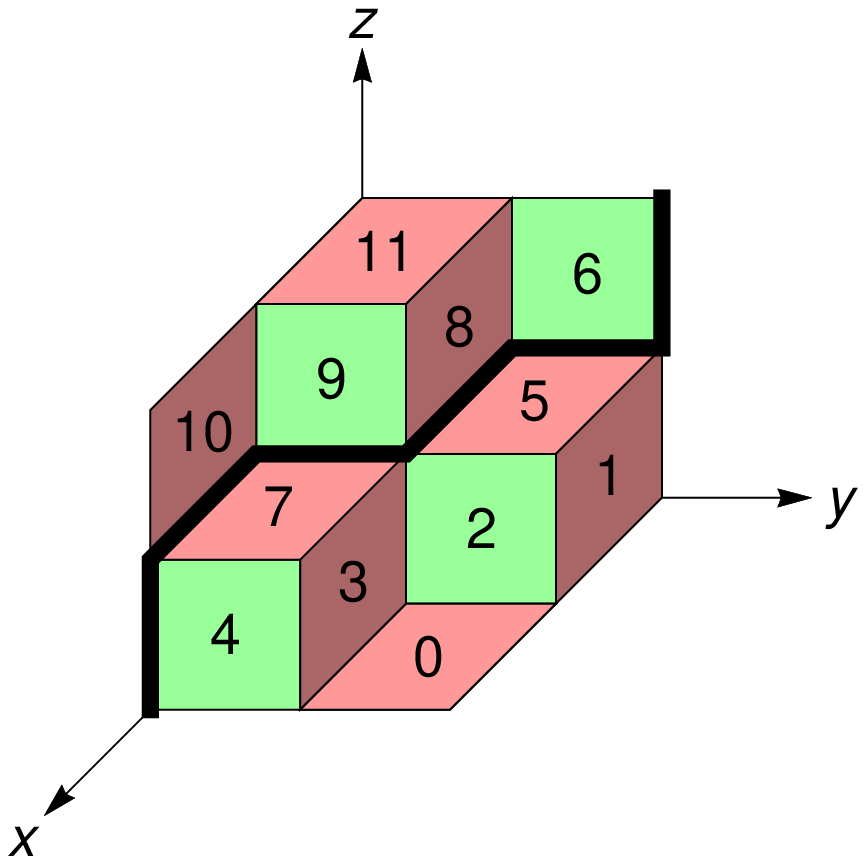}
\caption{\label{fig:orderb}}
\end{subfigure}
\caption{ {\bf Example labelings.}  a) The left labeling follows the reference ordering, \eref{refpart}. b) The right labeling results from one binary flip, which uniquely fixes the position of a part boundary, visually represented by a dark line, separating the tiles in the left and right parts of \eref{fixedPart}.}
\label{fig:order}
\end{figure}

The representation of an octagonal tiling as a hexagonal partition, itself built on a plane partition, provides for combinatorial shortcuts that avoid the total enumeration of octagonal tilings.
In this counting algorithm, we explicitly enumerate hexagonal tilings and solve the partition problems on them.

Given a hexagonal tiling, we first label the tiles such that the sequence of heights in that order of the labeling is non-decreasing.
An example of such an ordering can be seen in \fref{ordera}.
Next, that sequence of faces (tiles) can be partitioned such that successive parts have increasing heights.
For a hexagonal tilings of 12 tiles being lifted to height $d = 2$, such a partition could look like:
\begin{equation} \label{eqn:refpart}
 0~1~2~3~|~4~5~6~7~8~9~|~10~11
\end{equation}
where the vertical lines separate the tiles into parts, and are therefore called ``part boundaries.''
The tiles with labels in the first part are assigned $L = 0$, the second part $L = 1$, and the third $L = 2$.
The result is an octagonal tiling, as demonstrated in \fref{oct}.

If $d$ is the height of the partition, then the number of partitions is:
\begin{equation} \label{eqn:part}
P(K,d) = \left( K + d \atop d \right)
\end{equation}

However, the partial ordering that constrains the labeling admits multiple valid labelings of faces.
Another example labeling can be found in \fref{orderb}.
We can list the labels in the same tile order as \eref{refpart} and \fref{ordera}.
\begin{equation} \label{eqn:altlabel}
 0~1~2~3~4~5~7~6~8~9~10~11
\end{equation}
which differs from labeling \ref{eqn:refpart} by a binary exchange of adjacent labels.

To enumerate every possible hexagonal partition, each valid ordering must be counted with \eref{part}.
However, partitions of different sequences can result in the same hexagonal partition.
For example, the partitions
\begin{equation}
\begin{aligned}
 0~1~2~3~|~4~5~6~7~8~9~|~10~11 \\
 0~1~2~3~|~4~5~7~6~8~9~|~10~11
\end{aligned}
\end{equation}
of labelings \eref{refpart} and \eref{altlabel}, respectively, represent the same hexagonal partition, resulting in double counting.
To prevent double counting, a part boundary must be fixed in one of the labelings such that the partitions based on the two orderings are unique:
\begin{equation} \label{eqn:fixedPart}
\begin{aligned}
 0~1~2~3~4~5~6~~~7~8~9~10~11 \\
 0~1~2~3~4~5~7~|~6~8~9~10~11
\end{aligned}
\end{equation}
which corresponds to the dark line in \fref{orderb}.
In general, in every pair of orderings that differ by a binary flip of labels one of the two orderings must contain a fixed part boundary between the flipped labels.
The ordering selected to contain the fixed part boundary is the ordering in which the flipped labels descend, that is $v_j > v_{j+1}$ where $v_i$ is the label in the $i$th position and $j$ is the position of the flip.

Since one part boundary is fixed, the number of possible partitions decreases.
In general, for an ordering of an hexagonal tiling with $j$ fixed part boundaries, there are
$$ P(K,d,j) = \left( K + d - j \atop d -j \right) $$
partitions.
The cardinality of the ensemble of octagonal tilings is therefore:
\begin{equation} 
\Omega = \sum_{h \in \{3\rightarrow2\}} P(K,d,j_h) = \sum_{h \in \{3\rightarrow2\}} \left( K + d - j_h \atop d -j_h \right)
\end{equation}
where $\{3\rightarrow 2\}$ is the set of valid labeling of a hexagonal tilings, e.g. \eref{altlabel}, and $j_h$ is the number of fixed part boundaries in $h$.
\begin{df}[$a_j$]
Let $ a_j(a,b,c) $ be the number of labelings of $(a,b,c)$ hexagonal tilings with $j$ fixed part boundaries.
%\begin{prop}
Therefore,
\begin{equation} \label{eqn:W}
 \Omega = \sum_j a_j(a,b,c)  \left( K + d - j \atop d -j \right) = \sum_j a_j(a,b,c)  \left( K + d - j \atop K \right)
\end{equation}
\end{df} 
%\end{prop}
This is the expression found by Destainville {\it et al.}, but with the number of fixed part boundaries replacing the notion of descents~\cite{Dest01}.
Symmetry provides the relation $ a_j(a,b,c) = a_j(c,a,b) $.
The MacMahon formula, \pref{macmahon}, can be applied to the subsequence of zero fixed part boundaries, i.e. zero descent:
$$ a_0(a,b,c) = \prod_{i=1}^a \prod_{j=1}^b \prod_{k=1}^c \frac{i + j + k - 1}{i + j + k - 2} $$

\begin{df} \label{def:element-wise}
We separate the terms of \eref{W} to define element-wise cardinalities and entropy densities:
\begin{align}
 \omega_j(a,b,c,d) &= a_j(a,b,c) \left(ab + bc + ac + d - j \atop d - j \right) \\
 s_j(a,b,c,d) &= \frac{1}{N} \ln  \omega_j(a,b,c,d) 
\end{align}
such that
$$ \Omega(a,b,c,d) = \sum_{j = 0}^d \omega_j(a,b,c,d) $$
$$ \forall j, S(a,b,c,d) \ge s_j(a,b,c,d) $$
\end{df}

A general procedure by which to compute $j_h$ from the ordered labeling of the hexagonal tiling $h$ is given by the descent theorem~\cite{Dest01}; we will summarize it here.
When two labels are flipped, a fixed part boundary should be counted in only one of the corresponding labelings.
To keep track of which labeling is assigned the fixed part boundary, an auxiliary ordering is introduced such that one of the labelings satisfies the ordering and the other does not.
The labeling that does not satisfy the ordering is assigned the fixed part boundary.
The break in the ordering is referred to as a descent, hence the descent theorem.
The constraints on the auxiliary ordering and the proof of the method are provided by Destainville et al~\cite{Dest01}.
Of particular consequence is that there exists exactly one zero descent labeling per hexagonal tiling, so $a_0(a,b,c)$ count hexagonal tilings.

\section{Calculation}
The problem has been reduced to the enumeration of hexagonal tilings and the counting of the fixed part boundaries, that is descents, of their labellings.
Here, we present a method for counting the fixed part boundaries over all tilings.

We could build the hexagonal tilings tile by tile, labeling tiles in the order they were added, and count each time we add a tile that breaks the reference ordering, resulting in a fixed part boundary.
The construction follows the partial ordering of the tiles in the hexagonal tiling, and thus results in decisions both in which tile to add and where to add it.
If these decisions are represented as a directed tree, then the leaves will enumerate labeled hexagonal tilings and the internal vertices will enumerate labeled partial tilings.
%Whether or not adding a tile adds induces a fixed part boundary depends on the vertex representing the new partial tiling, it's parent, and its grand-parent.
We will continue to use the representation of the algorithm on a graph.

The partial ordering of the tiles in the hexagonal tiling can be thought of as a boundary on the partial tilings.
Additional tiles can be added so long as they fit within the external hexagonal boundary and share an edge with the internal partial tiling boundary.
Thus, the decision of which tile to add and where to add it depends only on the boundary of the partial tiling, and not on its interior.  
Only two pieces of information from the interior of the tiling are relevant: the number of fixed part boundaries thus far, and the last tile added.
This allows for multiple partial tilings to be coalesced into the same vertex, reducing the size of the problem dramatically.
Instead of a partial tiling, each vertex stores a boundary, a list of last tiles added, and, for each last tile, a list of how many partial tilings have produced $j$ fixed part boundaries.
Alternatively, the list of last tiles added can be considered the set of edges that terminate in the vertex and the list of partial tiling counts with $j$ fixed part boundaries can be associated with the edges.

This can be considered dynamic programming:
the irrelevance of the details of the interior of the partial tiles leads to overlapping subproblems characterized by the boundary of the partial tiling and the most recently added tile.
Unlike many standard dynamic programming solutions, the recurrence structure forms a more general graph than the standard n-table.
We will now define the recurrence relation and base case, and discuss properties of the graph induced by the recurrence.

To avoid confusion with fixed part boundaries, we will henceforth refer to partial tiling boundaries as \textit{paths}, which uniquely identify vertices in our graph.
Each path is a sequence of edges in the cubic lattice that originate at the origin and terminate at the point $(a,b,c)$.
We represent each path as a sequence of unit vectors, each denoting the next edge in the path.
For example, the path:
\begin{equation}
P_t = \left( \hat{e}_0, \dots, \hat{e}_0, \hat{e}_1, \dots, \hat{e}_1, \hat{e}_2, \dots, \hat{e}_2 \right) \equiv \left(0, \dots, 0, 1,\dots, 1, 2,\dots, 2\right)
\label{eqn:path}
\end{equation} 
is the root of the tree, runs from the origin to $(a,0,0)$ to $(a,b,0)$ to $(a,b,c)$, and corresponds to a partial tiling with zero tiles.
The operation of adding a tile to the path is a binary flip of two elements of the sequence, as can be seen in \fref{graph}.
This relation on the paths forms a directed acyclic graph that we will call a \textit{sort graph}.
\begin{df}[Sort graph]
The $a,b,c$ sort graph is a directed graph with vertices $v \in \mathbb{N}^{a+b+c}$ such that
\begin{itemize}
  \item $a = \left| \left\{i : v_i = 0\right\} \right|$
  \item $b = \left| \left\{i : v_i = 1\right\} \right|$
  \item $c = \left| \left\{i : v_i = 2\right\} \right|$
\end{itemize}
A directed edge goes from $v$ to $v'$ iff $\exists j$ such that
\begin{itemize}
  \item $v_j = v'_{j+1} \ne v_{j+1} = v'_j$ 
  \item $\forall i \notin \{j,j+1\}, v_i = v'_i$
  \item $\sum_i v_i 3^i < \sum_i v'_i 3^i$
\end{itemize}
The vertices correspond to boundaries of partial tilings and the edges correspond to added tiles.
\end{df}
\begin{figure}
\begin{center}
$$\xymatrix@C=0em{
&& \ar@/_/[dl]_\mu \left(0, 0, 1, 1, 2, 2\right) \ar@/^/[dr]^\nu&& \\
& \left(0, 1, 0, 1, 2, 2\right) \ar@/_/[dl] \ar[d] \ar@/_/[dr] & &  \ar@/^/[dl] \ar[d] \ar@/^/[dr] \left(0, 0, 1, 2, 1, 2\right) & \\
\left(1, 0, 0, 1, 2, 2\right) & \left(0,1,1,0,2,2\right) & \left(0, 1, 0, 2, 1, 2\right) & \left(0, 0, 1, 2, 2, 1\right) & \left(0, 0, 2, 1, 1, 2\right) 
}$$ 
\end{center}
\caption{ {\bf Example sort graph.} The edges $\mu,\nu$ correspond to the faces labeled $0,1$ in \fref{order}, respectively.}
\label{fig:graph}
\end{figure}
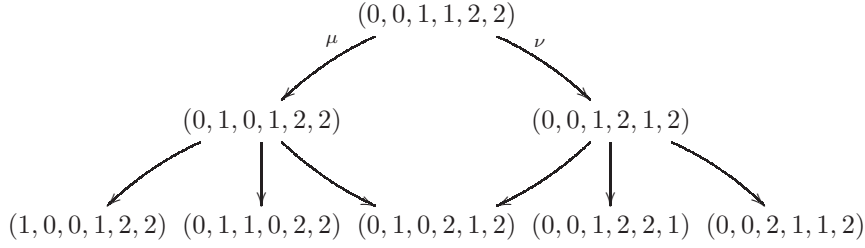
To count fixed part boundaries, we first give the faces of the cubic lattice an additional total ordering, $<$.
When the tile addition ordering breaks this total ordering, that is when tile $t_1$ is added after tile $t_2$ for $t_1 < t_2$, a descent is counted.
The total ordering is valid if each hexagonal tiling has exactly one zero-descent labeling.
We use the ordering of Destainville et al, which was used to prove the prior point.
\begin{df}[Reference ordering]
Let $t : E \rightarrow T$ be a mapping from edges in a sort graph to tiles of a hexagonal tiling, that is faces in a unit cubic lattice.
We represent the tiles by the coordinates of the point at their center.
Let $e = (u,v)$ and $j$ be the first index for which $u,v$ differ.
\begin{equation}
\vec{t}[e] = \frac{1}{2} (u_j + v_j) + \sum_i^{j-1} u_i
\end{equation}
where the elements of the sequences $u,v$ are taken to be unit vectors.

For two edges $e_1, e_2$:
\begin{equation}
e_1 < e_2 \text{ iff } \begin{cases} 
\vec{t}[e_1]_2 < \vec{t}[e_2 ]_2 & \text{ if } \neg (\vec{t}[e_1]_2 = \vec{t}[e_2 ]_2) \\
\vec{t}[e_1]_1 > \vec{t}[e_2 ]_1 & \text{ if } (\vec{t}[e_1]_2 = \vec{t}[e_2 ]_2) \wedge \neg (\vec{t}[e_1]_1 = \vec{t}[e_2 ]_1) \\
\vec{t}[e_1]_0 > \vec{t}[e_2 ]_0 & \text{ if } (\vec{t}[e_1]_2 = \vec{t}[e_2 ]_2) \wedge (\vec{t}[e_1]_1 = \vec{t}[e_2 ]_1) \wedge \vec{t}[e_1]_2 \in \mathbb{Z} \\
\vec{t}[e_1]_0 < \vec{t}[e_2 ]_0 & \text{ else } \\
\end{cases}
\end{equation}
\end{df}
Sort graphs have a layered structure.
Every path from the root to a vertex $v$ has the same length, which defines the depth of the vertex, $d(v)$.
There is a single vertex with zero out-degree.
The depth of this terminal vertex is the height of the sort graph, which is the number of tiles in the hexagonal tiling, $K$.

The process of counting the descents of the walks through the digraph is accomplished with a single breadth-first operation.
To each edge in the digraph, $e$, we associate a sequence $a^e_j$ which stores the number of walks with $j$ descents that are on that edge.
The first two edges, seen in \fref{graph}, are the base case:
\begin{equation}
a^\mu_j = a^\nu_j = \begin{cases} 1 & \text{ if } j = 0 \\ 0 & \text{ else } \end{cases}
\end{equation}
The edges add their sequences to their descendants, incrementing the index $j$ if the ordering on the edges is broken.
\begin{equation} \label{eqn:recurrence}
a^e_j = \sum_{e'} \begin{cases} a^{e'}_j & \text{ if } e \le e' \\ a^{e'}_{j-1} & \text{ else } \end{cases}
\end{equation}
where $e'$ runs over the edges with the destination vertex that is the source vertex of $e$.
Ultimately, the sequences of the two edges that terminate in the last path are summed, yielding the counts of labellings with $j$ fixed part boundaries: $a_j = \sum_t a^{t}_j$, where $t$ are the edges that terminate in the terminal vertex, that is the vertex with an out degree of zero.
To compute the cardinality of the ensemble, $\Omega(a,b,c,d)$, we need only $j \le d$.

Sort graphs are large.
The number of vertices is:
\begin{equation}\label{eqn:nodes}
V = \binom{a+b+c}{b+c} \binom{b+c}{c} \equiv \left( a+b+c \atop b ~ , ~ c \right)
\end{equation}
The number of edges is:
\begin{equation}\label{eqn:edges}
E = \frac{ab + ac + bc}{a+b+c} V
\end{equation}
The number of vertices in the $i$th layer of the graph are given by:
\begin{equation}
V_i = \sum_{j = 0}^i p(a,b+c;j) p(b,c;i-j)
\end{equation}
where $p(N,M;n)$ is the number of partitions of $n$ with no more than $N$ parts of size less than or equal to $M$.

The time complexity, $T(n)$, goes as the product of the number of edges and the payload size.
Let $n \sim a \sim b \sim c \sim d$.
The payload has a maximum size of $d$, so:
\begin{equation}
T(n) \sim n^2 \left( 3 n \atop n ~,~ n \right)
\end{equation}
Applying Stirling's approximation:
\begin{equation}
T(n) \approx \frac{n \sqrt{3}}{2 \pi} 3^{3n}
\end{equation}
The number of octagonal tiles, $N \sim n^2$.
Therefore, $T(N) \in O(\sqrt{N} \exp[\sqrt{N}])$.
Compared to the $O(\exp[N])$ complexity given for the algebraic formulation~\cite{Dest04}, this method is much more efficient.

The space complexity, $F(n)$, goes as the product of the payload and the size of the largest layer because the operation is breadth-first.
There are $K \sim n^2$ layers with the central layer being largest:
\begin{equation}
F(n) \sim n V_{\text{max}} = n \sum_{j=0}^{3 n^2/2} p(n,2n;3n^2/2-j) p(n,n;j)
\end{equation}

\subsection{Computational}

The calculation has been implemented in C/pthreads~\cite{Mueller1993} and Hadoop MapReduce~\cite{Bialecki2005}.
These implementations are available for use\footnote{\texttt{https://github.com/maxhutch/big-data-project}}.

The C/pthreads implementation follows a simple queue-based breadth-first operation supplemented with a hash table for searching.
The base-case vertices are added to a queue.
While the first element in the queue is not the terminal vertex, a vertex is dequeued and processed.
Processing consists of enumerating the binary flips in its path that correspond to outgoing edges.
For each outgoing edge, the destination vertex is initialized and added to the queue, if not there already, and the payload $a_j$ is updated following \eref{recurrence}.
After processing is complete, the vertex is deleted.
The threaded algorithm simply uses a thread-safe queue and hash table, adds a mutex to each vertex, and waits if the queue is empty.

The two essential features of this algorithm are a concise representation of a vertex and an efficient hash table.
The has function used is ternary interpretation of the identifying path, which is a perfect hash over the set of vertices.
However, only a subset comprised of no more than two layers of vertices are simultaneously resident in the hash table, so the performance is not guaranteed to be optimal.
Further, real memory systems provide strong incentives for locality, which is not accounted for with this hash.

The vertex consists of a identifying path, a list of preceding edges, and list of $a_j$ for each such edge.
The path is represented as an array of bit-pairs, which is $\log(3)/\log(4) \approx  80\%$ efficient.
The edges are identified by the coordinates of the center of the corresponding face in the cubic lattice.
The coordinates are integer or half-integer, so they are doubled and stored as a 3-tuple of bytes.
The $a_j$ values are not bounded and frequently overflow 64-bit integers.
The GNU multiprecision library is used to store them as arbitrary length integers, at the cost of fragmenting memory space.

The vast size and low edge-density of sort-graphs make them prime candidates for the MapReduce paradigm.
The MapReduce implementation operates on edges instead of vertices.
The \textit{map} procedure takes an edge $e' = (u,v)$ to a set of edges $\{e = (v,w)\}$, each of which contain the contribution of $a^{i'}$ to the summation for $a^e$, \eref{recurrence}.
The \textit{reduce} procedure sums the $a_j$ sequence over identical edges.
The map-reduce pair is repeated for each layer of the sort graph, implying a barrier between them.
At each iteration the edges are sorted by their destination vertex using the Hadoop \textit{partition}.
The sorting partitioner causes the majority of edges identified by the same path to be produced by the same process.
This locality is taken advantage of with the Hadoop \textit{combine} operation, which is identical to reduce but performed only on the subset of computer-local edges.
The use of partition and combine isolate the non-locality of the implementation to a single global sort.

\section{Results}

\subsection{Enumerative}

\begin{table}[h]
  \centering
    \begin{tabular}{c | c  c  c  c } 
$\alpha$ & $\beta = 1$ & $\beta = 2$ & $\beta = 3$ & $\beta = 4$ \\
\hline
1 & 0.34657359 & 0.33313333 & 0.30562875 & 0.27921012 \\
2 & 0.35795838 & 0.33975183 & 0.31110489 & 0.28371882 \\
3 & 0.35978800 & 0.34128074 & 0.31260400 & 0.28493173 \\
4 & 0.36021813 & 0.34192781 & 0.31325143 & 0.28543924 \\
5 & 0.36031464 & 0.34227989 &            &            \\
6 & 0.36031361 & 0.34249774 &            &            \\
7 & 0.36028282 &            &            &            \\
8 & 0.36024459 &            &            &            \\
    \end{tabular}

  \caption{ {\bf Entropy densities.} Entropy density for $\beta = 1, 2, 3, 4$.  Blanks correspond to numerically intractable cases.  Entropy density is taken to be the logarithm of the number of tilings over the number of tiles: $\ln \Omega / N$.}
  \label{tbl:sumres}
\end{table}

We have been able to expand on the numerical results of Destainville {\it et al.}.
\tref{sumres} provides entropy densities, \eref{S}, for the first four integral $\beta$ cases.
The $\beta = 1$ column corresponds to the `diagonal' case.
Surprisingly, it is not a monotonically convergent sequence, as it decreases slowly after an initial rapid increase.
The full integer sequence $\Omega(\alpha,1)$ can be found on the Online Encyclopedia of Integer Sequences\footnote{\texttt{http://oeis.org/A093937}}.

\begin{table}[h]
{\scriptsize
  \begin{tabular}{ c | r  r  r  r  r  r r  r}
    j & $\alpha = 1$ & $\alpha = 2$ & $\alpha = 3$ & $\alpha = 4$ & $\alpha = 5$ & $\alpha = 6$ & $\alpha = 7$ & $\alpha = 8$ \\
 \hline
0   &   0.346667    &   0.312917    &   0.281481    &   0.259063    &   0.242667    &   0.230231    &   0.220442    &   0.212526    \\
1   &   0.000000    &   0.331667    &   0.327778    &   0.303333    &   0.282400    &   0.265833    &   0.252653    &   0.241875    \\
2   &       &   0.273333    &   0.346296    &   0.332604    &   0.312267    &   0.294306    &   0.279252    &   0.266719    \\
3   &       &       &   0.338889    &   0.349896    &   0.334800    &   0.317546    &   0.301905    &   0.288385    \\
4   &       &       &       &   0.354063    &   0.350533    &   0.336296    &   0.321224    &   0.307422    \\
5   &       &       &       &       &   0.358333    &   0.350648    &   0.337551    &   0.324141    \\
6   &       &       &       &       &       &   0.359676    &   0.350714    &   0.338672    \\
7   &       &       &       &       &       &       &   0.360068    &   0.350885    \\
8   &       &       &       &       &       &       &       &   0.360182    \\
  \end{tabular}

}
  \caption{{\bf Element-wise entropies.} $s_j(\alpha,1)$ as in \dref{element-wise}.}
  \label{tbl:sj}
\end{table}

We also present the element-wise entropies, $s_j(\alpha, 1)$, in \tref{sj}.
The complete sequences $a_j(\alpha,1)$ can be found on the Online Encyclopedia of Integer Sequences \footnote{\texttt{http://oeis.org/A217311}}.
The entropy is the logarithm of a sum, which can be approximated as the logarithm of the maximal element of the sum if the number of terms is small.
For $\alpha \le 3$, the maximal element corresponds to $j = \alpha - 1$.
For $\alpha > 3$, the maximal element shifts to $j = \alpha$ and the dominance increases as $\alpha$ increases.

\subsection{Computational}

The low level implementation requires less aggregate time than MapReduce for computing $\Omega$ for $\alpha \le 8$.
The $\beta=1, \alpha=8$ case required nearly 14 days at over 1TB of shared memory on the Blacklight supercomputer at the Pittsburgh Supercomputing Center\footnote{\texttt{http://www.psc.edu/index.php/computing-resources/blacklight}}.
However, the MapReduce implementation is very effective when computing the entire $a_j$ sequence, not just $j \le \alpha$.
This is because MapReduce is latency bound in our case, so the marginal cost of increasing the computational workload and message size is low.
Given the interesting properties of the $a_j$ sequences, independent of their role in computing $\Omega$, MapReduce is a promising tool for future work.

%The $\alpha = 9$ case has $325/12\approx 27$ times more edges than $\alpha=8$, and the runtime and memory usage are $\Omega(E)$, so this computation is likely to be intractable for the next few years.

\section{Discussion}
\begin{figure}
\resizebox{\columnwidth}{!}{\input{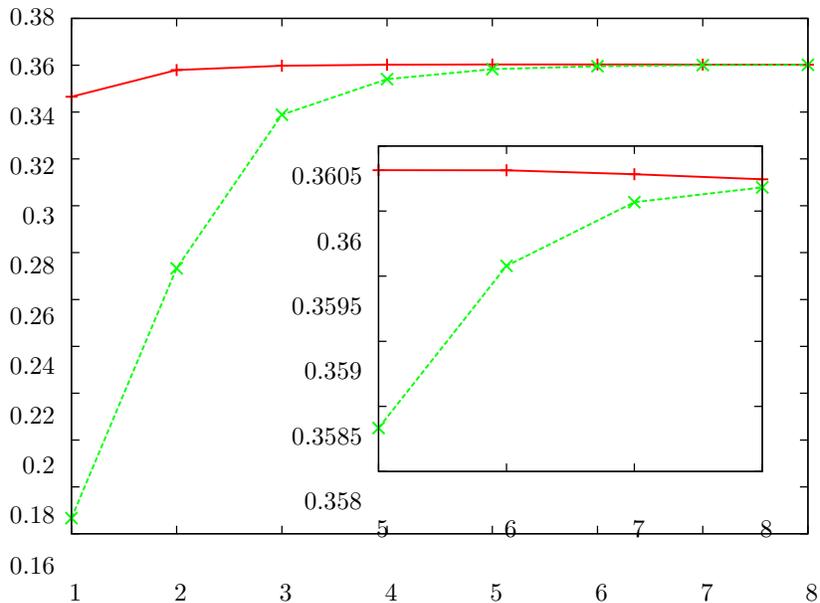}}
\caption{ {\bf Convergence wrt $\alpha$.} 
Convergence of the entropy density $S(\alpha,1)$, solid, and lower bound $s_\alpha(\alpha,1)$, dashed, with respect to the problem size $\alpha$.
Inset focuses on the break in monotonicity of $S(\alpha,1)$ seen at $\alpha > 5$.
The entropy density is the first column of \tref{sumres} and the lower bound is the sub-diagonal of \tref{sj}.}
\label{fig:S}
\end{figure}

The results suggest limiting behavior of the sequences $a_j$ and $s_j$.
First, the diagonal of \tref{sj} appears to be monotonic increasing.
Because every element of the columns of \tref{sj} are lower bounds on the entropy density $S(\alpha,1)$, this establishes a monotonically increasing lower bound.
\begin{hyp} \label{hyp:lower}
$$ \forall \alpha, s_{\alpha}(\alpha, 1) < s_{\alpha+1}(\alpha+1, 1) $$
or, equivalently: 
\begin{equation*}
\forall \alpha, \left(a_\alpha(\alpha, \alpha, \alpha)\right)^{(\alpha+1)^2 / \alpha^2} < a_{\alpha+1}(\alpha+1, \alpha+1,\alpha+1) 
\end{equation*}
\end{hyp}

A weaker suggestion is that beyond the peak at $\alpha = 5$, the entropy density is monotonically decreasing.
This establishes a monotonically decreasing upper bound on the entropy density in the thermodynamic limit, $S(\infty,1)$.
\begin{hyp} \label{hyp:upper}
$$ \forall \alpha > 6, S(\alpha+1,1) < S(\alpha,1) $$
\end{hyp}

Taking the bounds together, one can estimate $S(\infty,1)$ to 5 digits.
The convergence of these two bounds can be seen in \fref{S}.
\begin{prop}[Bounds]
If \yref{lower} and \yref{upper}, then we can bound the thermodynamic limit above by $S(8,1) = 0.36024459\ldots$ and below by $s_8(8,1) = 0.360182\ldots$, resulting in $S(\infty,1) = 0.36021 \pm 0.00003$.
\end{prop}

Finally, after some finite size effects at $\alpha \le 3$, the columns of \tref{sj} appear to be monotonically increasing.
\begin{hyp} \label{hyp:ddom}
$$\forall \alpha > 3, x_1 < x_2 \le \alpha, s_{x_2}(\alpha, 1) > s_{x_1}(\alpha, 1) $$
\end{hyp}

\section{Conclusions}
The algorithm of Destainville {\it et al.} can be stated in terms of fixed part boundaries identified by descents in directed labelings over the ensemble of hexagonal tilings.
The complexity of the resulting algorithm is $O(\exp[\sqrt{N}])$, where $N$ is the number of octagonal tiles, which compares favorably to the $O(\exp[N])$ complexity of recent alternatives.

Results for octagonal ensembles, $\Omega(a,b,c,d)$, are constructed from a rich space of sequences $a_j(a,b,c)$.
Enumerative results through $\alpha = 8$ strongly suggest relationships in the $a_j(a,b,c)$ sequences that can be used to construct bounds for the thermodynamic limit of the entropy density $S(\infty,1)$.
Assuming \yref{lower} and \yref{upper}, one can estimate the thermodynamic limit as $S(\infty,1) = 0.36021 \pm 0.00003$.

Proofs of \yref{lower}, \yref{upper}, and \yref{ddom} are open problems.
The formula for the entropy density allows for significant approximation of the $a_j(a,b,c)$ sequences while still converging to the correct thermodynamic limit.
Stepping stones towards such an approximation include expressions for $\sum_j a_j(a,b,c)$ or lower bound for $(a_\alpha(\alpha,\alpha,\alpha))^{(\alpha+1)^2/\alpha^2}$.

\section*{Acknowledgements}
We would like to acknowledge useful conversations with Rhys Povey and Jonathan Paulson.

\bibliographystyle{elsarticle-num} 
\bibliography{../bib/library}

\end{document}